\documentclass[12pt]{article}
\usepackage{amsfonts,amssymb,amsmath,latexsym}
\usepackage{epsfig,graphicx}
\newtheorem{theorem}{Theorem}
\newtheorem{lemma}[theorem]{Lemma}

\newcommand{\pf}{\noindent \mbox{{\bf Proof}: }}
\begin{document}
\title{Coincidence of  length spectra does not imply
isospectrality}
\author{S. A. Fulling and Peter Kuchment\\
Mathematics Department, Texas A\&M University\\ College Station,
TX 77843-3368} \maketitle

\begin{abstract}
Penrose--Lifshits mushrooms are planar domains coming in
nonisometric pairs with the same geodesic length spectrum. Recently
S.~Zelditch raised the question whether such billiards also have the
same eigenvalue spectrum for the Dirichlet Laplacian (conjecturing
``no''). Here we show that generically (in the class of smooth
domains) the two members of a mushroom pair have different spectra.
\end{abstract}

\section{Introduction}

 Michael Lifshits (unpublished), exploiting a type of construction
attributed to R.~Penrose (see, e.g., \cite{Rau}), constructed a
class of pairs of planar domains that, while not isometric, have
periodic geodesics of exactly the same lengths (including
multiplicities).
 At least when the boundaries are smooth ($C^\infty$),
 it follows that the two billiards have the same wave invariants,
in the sense that the traces of their wave groups,
$\,\cos(t\sqrt{\Delta})\,$,
 differ at most by a smooth function~\cite{Me}.
 In a recent review of the inverse spectral problem \cite{Zeld}
 S.~Zelditch asked whether the Dirichlet Laplacians, $\Delta$, for the
two domains are necessarily isospectral, judging that proposition
 ``dubious'' but not yet refuted.
 Given the refutation,
 such billiards  provide a kind of converse to the famous
examples of ``drums that sound the same'' \cite{GWW},
 being drums that sound different but are very similar
geometrically --- in fact, in the geometrical features deemed most
relevant to  spectrum.

 In this paper we show how to construct smooth Penrose--Lifshits mushroom
 pairs that are not isospectral, and we argue that inequality of
the Dirichlet spectra is, in fact, quite generic.
 Since the domains are smooth (but not convex), the spectral difference
 is not attributable to diffraction from corners, which would muddy
the definition or the relevance of ``periodic
geodesics''\negthinspace.

\section{Main result}

 The construction of a mushroom  starts from
 a half-ellipse $E$ with foci $F$ and~$\tilde F\,$:
\begin{center}
\includegraphics[width=4in,height=2in]{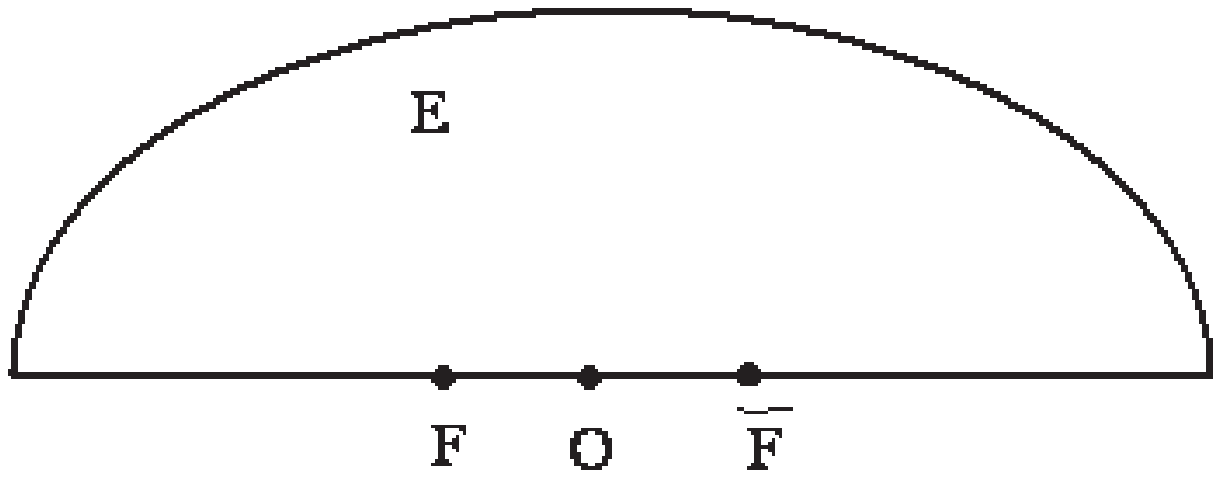}
\end{center}
 We use the tilde, whether applied to regions, curves, or points,
to indicate the operation of reflection through the minor axis of
the ellipse.
 If two entities are interchanged by that reflection, we call them
\emph{dual}.
Next, add two bumps, $B_1$ on the left and  $B_2$ on the
right, with $\tilde B_1 \ne B_2\,$,
 to form a smooth domain~$\Omega$:
\begin{center}
\includegraphics[width=4in,height=2in]{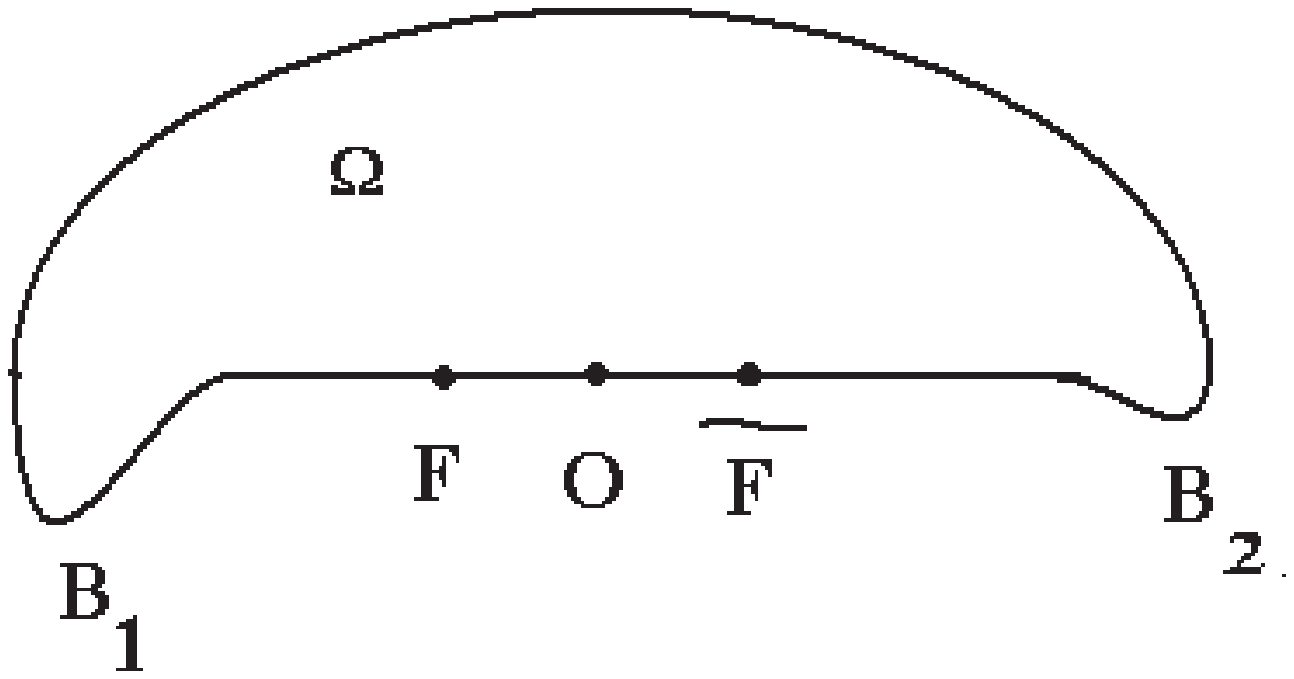}
\end{center}
 Finally, add another bump (not self-dual)
 between the foci in two dual
 ways ($M$ and $\tilde M$) to get two domains
$\Omega_j\,$:
\begin{center}
\includegraphics[width=4in,height=2in]{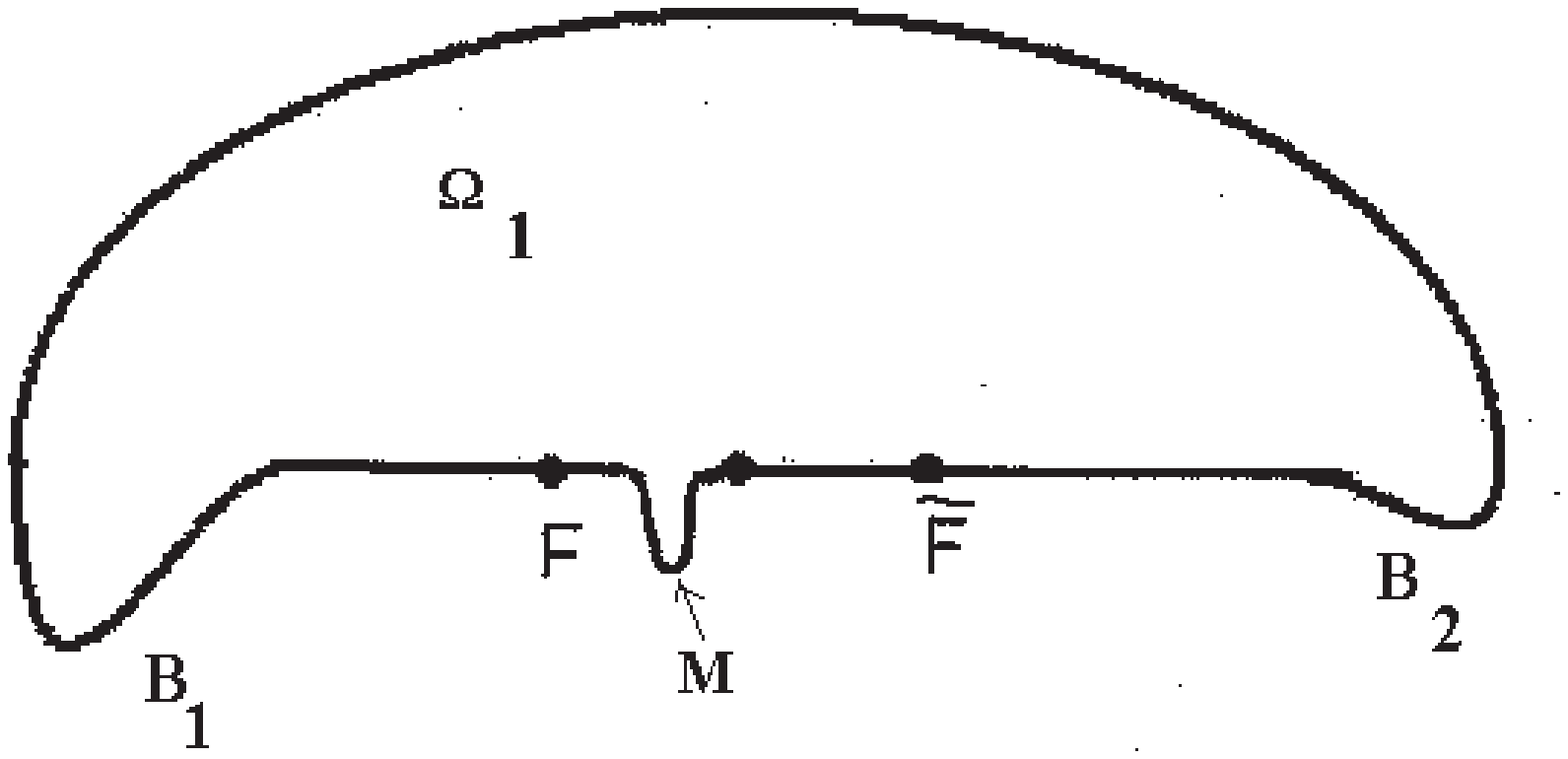}
\end{center}
\begin{center}
\includegraphics[width=4in,height=2in]{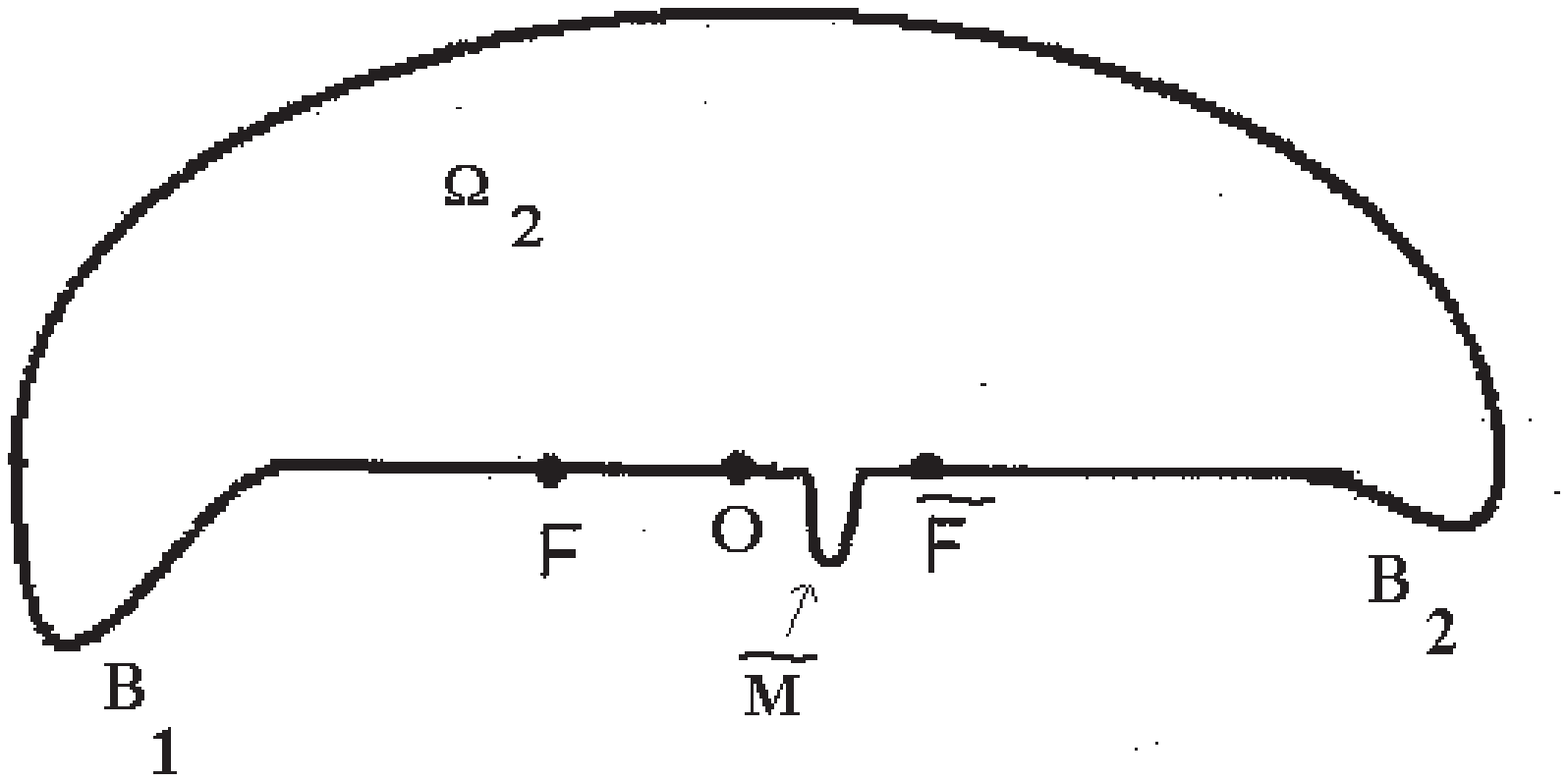}
\end{center}
We  call the domains $\Omega_1$ and $\Omega_2$
  constructed in this manner
\emph{a  Penrose--Lifshits mushroom pair}.

We repeat that the bumps can be added in such a way that the
boundaries remain smooth.
 That assumption, however, is needed mainly to draw a clean
conclusion about equality of the length spectra.
 The conclusions about the Dirichlet spectra hold even if the
domain has corners (in which case bump $B_2$ is superfluous).

\begin{theorem}\label{T:non-isospectral}
If $B_1$ and $B_2$ are given and not dual,
%mirror-symmetric with respect to the vertical axis of the ellipse,
 then there exist dual bumps
$M$ and $\tilde M$ such that the resulting Penrose--Lifshits
mushrooms
  $\Omega_j$
have the same length spectra and wave invariants but are not
isospectral.
\end{theorem}

\goodbreak\pf

First we review the proof that the length spectra coincide
 \cite{Me,Zeld}.
% of all, as it has already been mentioned, it is known that
%domains $\Omega_j$ have same lengths spectra and wave invariants
%\cite{Me,Rau,Zeld}. Although it is harder to prove the claim about
%the wave invariants, the argument about the geodesic lengths
%spectra is simple, and so we repeat it here. Indeed,
The geodesics in an ellipse fall into two disjoint categories
 \cite{KR,Rau,Berry}: those that intersect
the major axis between the foci, and those that do so at or beyond
the foci. (The only exception is the major axis itself.
The smoothness assumption guarantees that the major axis will not bifurcate
 in $\Omega_j$ by diffraction.)
 It follows that a similar division holds  for the domains
$\Omega_j$ we have just described:
any geodesic originating in a bump $B_1$ or $B_2$ can never reach a
bump $M$ or $\tilde M$, and vice versa.
  Now, the geodesics that do not
intersect the focal segment $F\tilde F$ are exactly the same for
the two
domains. On the other hand, those for $\Omega_1$ that do intersect
this segment are identified one-to-one
 with their duals in $\Omega_2$ by the
 reflection operation.
%symmetry with respect to the minor axis of the ellipse (we remind
%the reader that the ``tongues'' $M_j$ are symmetric with respect
%to that axis).
 This shows length isospectrality.
 Equality of the wave traces modulo smooth functions follows from
the hyperbolic propagation of singularities along geodesics ---
 see \cite{Me} and references therein.

 Our main task is  to show nonisospectrality for some choice
of  $M$. Consider the spectrum of $\Omega_1$
assuming that the bump $M$ is small
 and has support on the left half of the focal segment;
   i.e., to construct $\Omega_1$
  the (open) segment $FO$
 in the boundary of $\Omega$ is
perturbed by the graph of a smooth, compactly
supported (and nonpositive)  function $\epsilon f(x)$, where
$\epsilon$ is a small parameter.
  Let  $\psi_0$ be the ground
state of the Dirichlet Laplacian on $\Omega$ and $\lambda_0$ be
the corresponding lowest eigenvalue.
 The known Rayleigh--Hadamard formula
for change of the spectrum under domain perturbations (e.g.,
\cite{GS,Krein,Zeld} or \cite[Section 15.1, Exercise 9]{Garab})
 shows that if
%the shape perturbation function $f(x)$ is not orthogonal to the
%normal derivative of $\psi$ on the boundary of $\Omega$, i.e., if
\begin{equation*}
 \int\limits_{\partial \Omega} \left(\dfrac{\partial\psi_0
(x)}{\partial \nu}\right)^2 f(x) d\sigma(x) \neq 0,
 \end{equation*}
  where
$\dfrac{\partial\psi_0 (x)}{\partial \nu}$ is the normal
derivative of the eigenfunction on the boundary, then the lowest
eigenvalue $\lambda_0$ changes under the perturbation. In fact,
this integral gives the derivative at $\epsilon =0$ of the lowest
eigenvalue with respect to~$\epsilon$. Thus, if we could guarantee
that the values of this integral are different for the two small
perturbation domains $\Omega_j\,$, this would imply their
non-isospectrality: the lowest Dirichlet eigenvalues would
change with different velocities. Since the choice of the
perturbation shape $f$ is in our hands, in order to make these
integrals different, and thus domains nonisospectral, it is
sufficient to have two mutually dual segments
  inside the focal
segment such that the square of the normal derivative
 of the ground state  is not an even
 function on their union,~$I$.
 Indeed, in this case we could find an even perturbation $f$
that would provide nonequal integrals (in fact, almost any
perturbation would do).

\begin{center}
\includegraphics[width=4in,height=2in]{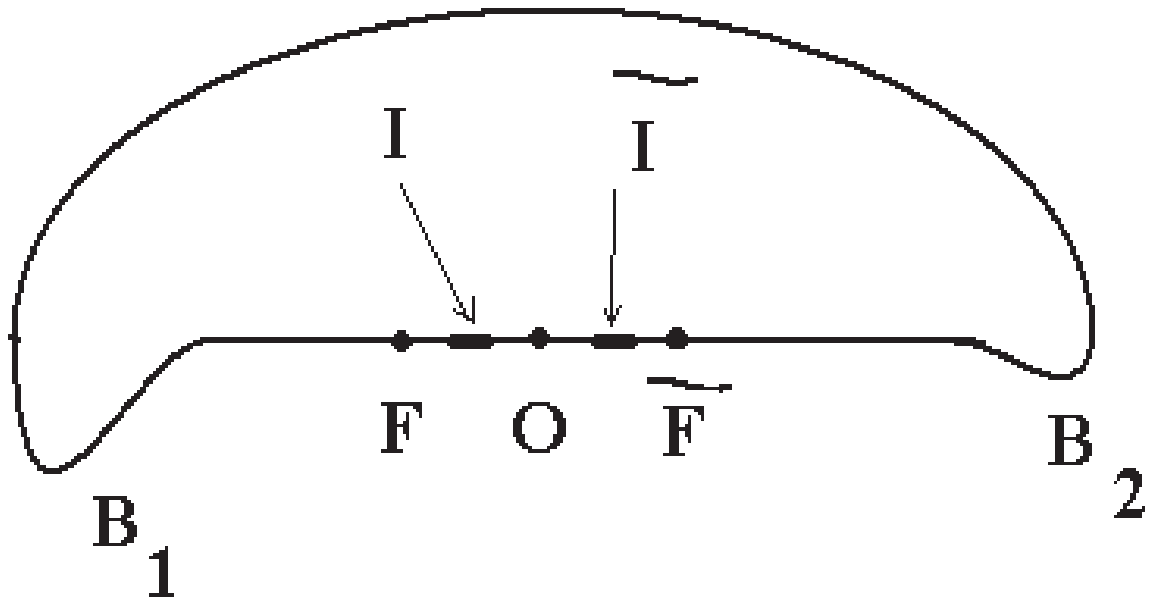}
\end{center}

\begin{lemma}\label{L:lemma}
If the bumps $B_j$ are not dual  with respect to the minor axis of
the ellipse, there is no self-dual union $I$ of two segments
inside $F\tilde F$ such that the square of the normal derivative
$\dfrac{\partial \psi_0}{\partial \nu}$ of the ground state
$\psi_0$ for $\Omega$ is even on~$I$.
\end{lemma}
{\bf Proof of the lemma.} Suppose that $\left(\dfrac{\partial
\psi_0}{\partial \nu}\right)^2$ is even on $I$. Since the normal
derivative is continuous, by shrinking $I$ if necessary,
%symmetric intervals can be made such that the normal derivative is
 we may assume that
 $\dfrac{\partial\psi_0}{\partial \nu}$
itself  is either even or odd on $I$.
Suppose first that the normal derivative is even. Introduce the
orthogonal cartesian coordinates centered at $O$
 %the ellipse's center
and with $x$-axis going along the major axis. Consider the
function $\psi_1=\psi_0(-x,y)$. Both $\psi_0$ and $\psi_1$
satisfy the same eigenfunction equation inside the half-ellipse
$E$ and have the same Cauchy data on~$I$. Therefore, according to
Holmgren's uniqueness theorem, they agree on their common domain.
In particular, $\psi_0$ must satisfy zero Dirichlet boundary
conditions not only on $\partial \Omega$, but also on its mirror
reflection with respect to the minor axis of the ellipse. Since
the bumps $B_j$ are assumed  not dual to each other, we conclude that
$\psi_0$ vanishes somewhere inside~$\Omega$
 (or $\psi_1$ somewhere inside~$\tilde\Omega$).
 That is, $\psi_0$ has a nodal curve, which is well known
to be impossible for a ground state (e.g., \cite{Courant,Garab}).
%
%An analogous consideration leads to contradiction
 If $\dfrac{\partial \psi_0}{\partial \nu}$ is odd on $I$, one only
needs to define $\psi_1$ as $-\psi_0(-x,y)$ to obtain an analogous
 contradiction.

This concludes the proof of the Lemma, and hence of  Theorem
\ref{T:non-isospectral}.
 \goodbreak

In fact, a closer look at the proof of the theorem shows
that the nonisospectrality holds for smooth Penrose--Lifshits
mushrooms
 $\Omega_j$ for any nondual bumps $B_j$ and for
``generic'' dual bumps $M, \tilde M\,$:

\begin{theorem}\label{T:generic}
For any fixed choice of nondual bumps $B_j\,$,
nonisospectrality holds for an open and dense (in
$C^\infty$-topology) set of Penrose--Lifshits
 pairs~$\Omega_j\,$.
\end{theorem}
\pf Indeed, the set of nonisospectral pairs $\Omega_j$ is
obviously open. The previous theorem states that the closure of
this set contains the domain $\Omega$ (i.e., the one where the
bumps $M, \tilde M$ are absent). To show density, one can apply a
similar proof by small perturbation to any pair of mushroom
domains $\Omega_j$ of the type constructed above. Indeed, if the
pair is already non-isospectral, there is nothing to prove. If it
is isospectral, let $\psi^j$ be the ground state in $\Omega_j$. As
in Theorem \ref{T:non-isospectral}, the perturbation method
described above works if one can show absence of a dual pair $J,
\tilde J$ of pieces of the boundaries $\partial \Omega_j$ such
that $J\subset
\partial M$, $\tilde J \subset \partial \tilde M$ and that
$\left(\dfrac{\partial \psi^1}{\partial \nu}\right)^2|_J$ is equal
after reflection to $\left(\dfrac{\partial \psi^2}{\partial
\nu}\right)^2|_{\widetilde J}$. Now, the same consideration as in
the proof of Lemma \ref{L:lemma} applies to justify this claim.

\section{Comments and acknowledgments}
\begin{itemize}
\item A different proof of generic non-isospectrality claimed
    in Theorem \ref{T:generic} follows from existence of
    non-isospectral mushroom domains (Theorem
    \ref{T:non-isospectral}), analytic dependence of the ground
    state on the domain \cite{GS}, and connectedness of the
    manifold of these domains.
\item As it is not hard to establish, the set of non-isospectral
mushroom pairs is open in a much weaker
    topology than $C^\infty$. Indeed, if the domains $\Omega_j$ are
    distorted by a pair of dual (in the sense used in this text)
    $C^2$-diffeomorphisms that are $C^2$-close to identity, the
    non-isospectrality is preserved.
\item One can find discussion of the effects of domain variation
for general elliptic boundary value problems in the nice little
book \cite{Krein}, which regrettably is available only in Russian.
Some of its results can be found in preceding publications of the
authors of that book.
\end{itemize}

This small note is the result of discussion at a working seminar
of the recent survey \cite{Zeld} by Steve Zelditch. The authors
would like to thank the seminar participants G.~Berkolaiko,
J.~Harrison and B.~Winn, as well as S.~Zelditch and J. Zhou,
 for discussion and useful comments.

The work of P.~Kuchment was partially supported by the NSF Grant
DMS 0406022. P.~Kuchment expresses his gratitude to NSF for this
support. The content of this paper does not necessarily reflect
the position or the policy of the federal government of the USA,
and no official endorsement should be inferred.

\end{document}